\documentclass[11pt,twoside]{article}
\usepackage{amsmath}
\usepackage{amsthm}
\usepackage{amssymb}
\usepackage{array}
\usepackage{geometry}
\usepackage{graphicx}
\usepackage{enumerate}
\usepackage{lscape,graphicx}
\usepackage{fancyhdr}
\usepackage{makeidx}
\usepackage{amsfonts}
\usepackage{graphicx}
\usepackage{CJK}
\usepackage{amsmath}
\usepackage{latexsym}
\usepackage{amssymb}
\usepackage{mathrsfs}
\usepackage{amsthm}
\usepackage{array}
\usepackage{enumerate}
\usepackage{longtable}
\usepackage{multirow}
\makeindex
\usepackage{wasysym}
\usepackage{chapterbib}

\normalsize
\begin{document}


\newpage
\setcounter{equation}{0}
\begin{center}
\vskip1cm
{\Large
\textbf{The Matrix-variate Dirichlet Averages and Its Applications} }
\vskip.5cm
%

\vspace{0.50cm}
\textbf{\large Princy T$^1$ and  Nicy Sebastian$^2$}\\

$^1$Department of Statistics, Cochin University of Science and Technology, Cochin\\
Kerala 682 022, India\\ {Email: {\tt{princyt@cusat.ac.in}}}\\

$^2$Department of Statistics,
St Thomas College, Thrissur, Kerala 680 001, India\\ {Email: {\tt{nicyseb@yahoo.com}}}\\

\end{center}

\begin{center} {\bf }  \vskip 0.10truecm \noindent \\

\vskip.5cm\noindent \centerline{\bf Abstract}
\end{center}
\par
This paper is about Dirichlet averages in the matrix-variate case or averages of functions over the Dirichlet measure in the complex domain. The classical power mean  contains the harmonic mean, arithmetic mean and geometric mean (Hardy, Littlewood and Polya), which is generalized to $y$-mean by deFinetti and hypergeometric mean by Carlson, see the references herein. Carlson's hypergeometric mean is to average a scalar function over a real scalar variable type-1 Dirichlet measure and this in the current literature is known as Dirichlet average of that function. The idea is examined when there is a type-1 or type-2 Dirichlet density in the complex domain. Averages of several functions are computed in such Dirichlet densities in the complex domain.  Dirichlet measures are defined when the matrices are Hermitian positive definite. Some applications are also discussed.
		 \\



{\section{Introduction}}
\thispagestyle{empty}

In Hardy et al.(1952) there is a discussion of the classical power mean, which contains the harmonic, arithmetic and geometric means. The classical weighted average is of the following form:
$$f(b)=[w_1z_1^b+...+w_nz_n^b]^{\frac{1}{b}}.
$$where all the quantities are real scalar where $w'=(w_1,...,w_n), z'=(z_1,...,z_n)$, $z_j>0,w_j>0,j=1,...,n$, $\sum_{j=1}^nw_j=1$ with a prime denoting the transpose. For $b=1$, $f(1)$ gives $\sum_{j=1}^nw_jz_j$ or the arithmetic mean; when $b=-1$, $f(-1)$ provides $[\sum_j(\frac{w_j}{z_j})]^{-1}=$ the harmonic mean and when $b\to 0_{+}$ then $f(0_{+})$ yields $\prod_{j=1}^nz_j^{w_j}=$ the geometric mean. This weighted mean $f(b)$ is generalized to $y$-mean by deFinetti [ deFinetti (1974)] and to hypergeometric mean by Carlson [Carlson (1977)]. A real scalar variable type-1 Dirichlet measure is involved for the weights $(w_1,...,w_{n-1})$ in Carlson's generalization, and then average of a given function is taken over this Dirichlet measure. In the current literature this is known as Dirichlet average of that function, the function need not reduce to the classical arithmetic, harmonic and geometric means.\\

The paper is organized as follows: Section $1.1$ gives the basic concepts for developing the theory of the matrix-variate Dirichlet measure in complex domain. Dirichlet averages for a function of matrix argument in the complex domain is developed in section $2$. In section $3$, we discuss the complex matrix-variate type-2 Dirichlet measure and averages over some useful matrix-variate functions. Rectangular matrix-variate Dirichlet measure is presented in section $4$. Some of the useful areas of applications are listed in section 5.

{\subsection{Complex Domain}}

\vskip.2cm In the present paper, we consider Dirichlet averages of various functions over  Dirichlet measures in the complex domain in the matrix-variate cases. All matrices appearing in this paper are Hermitian positive definite and $p\times p$ unless stated otherwise. In order to distinguish,  matrices in the complex domain will be denoted by a tilde as $\tilde{X}$ and real matrices will be written without the tilde as $X$. We consider real-valued scalar functions of the complex matrix argument and such functions will be averaged over complex matrix-variate Dirichlet measure. The following standard notations will be used: ${\rm det}(\tilde{X})$ will mean the determinant of the complex matrix variable $\tilde{X}$. The absolute value of the determinant will be denoted by $|{\rm det}(\cdot)|$. This means that if ${\rm det}(\tilde{X})=a+ib,i=\sqrt{-1}$ then $\sqrt{(a+ib)(a-ib)}=(a^2+b^2)^{\frac{1}{2}}=|{\rm det}(\tilde{X})|$. ${tr}(\cdot)$ will denote the trace of $(\cdot)$. $\int_{\tilde{X}}$ is integral over all $\tilde{X}$ where $\tilde{X}$ may be rectangular, square or positive definite. $\tilde{X}>O$ means that the $p\times p$ matrix $\tilde{X}$ is Hermitian positive definite. Constant matrices, whether real or in the complex domain will be written without the tilde unless the fact is to be stressed. and in that case we use a tilde. $O<A<\tilde{X}<B$ means $A>O,\tilde{X}-A>O,B-\tilde{X}>O$ where $A$ and $B$ are $p\times p$ constant positive definite matrices. Then $\int_{O<A<\tilde{X}<B}f(\tilde{X}){\rm d}\tilde{X}=\int_A^Bf(\tilde{X}){\rm d}\tilde{X}$ means the integral over the Hermitian positive definite matrix $\tilde{X}>O$ such that $O<A<\tilde{X}<B$ and $f(\tilde{X})$ is a real-valued scalar function of matrix argument $\tilde{X}$ and ${\rm d}\tilde{X}$ stands for the wedge product of differentials, that is, for $\tilde{Z}=(\tilde{z}_{ij})=X+iY$, a $m\times n$ matrix of distinct variables $\tilde{z}_{ij}$'s, where $X$ and $Y$ are real matrices, $i=+\sqrt{-1}$, then the differential element ${\rm d}\tilde{Z}={\rm d}X\wedge{\rm d}Y$ with ${\rm d}X$ and ${\rm d}Y$ being the wedge product of differentials in $X$ and $Y$ respectively. For example, ${\rm d}X=\wedge_{i=1}^m\wedge_{j=1}^n{\rm d}x_{ij}$, if $X=(x_{ij})$ and $m\times n$. When $\tilde{Z}$ is Hermitian then $X=X'$ and $Y=-Y'$. In this case ${\rm d}X=\wedge_{i\ge j=1}^p{\rm d}x_{ij}=\wedge_{i\le j=1}^p{\rm d}x_{ij}$ and ${\rm d}Y=\wedge_{i<j=1}^p{\rm d}y_{ij}=\wedge_{i>j=1}^p{\rm d}y_{ij}$. The complex matrix-variate gamma function will be denoted by $\tilde{\Gamma}_p(\alpha)$, which has the following expression and integral representation:

$$\tilde{\Gamma}_p(\alpha)=\pi^{\frac{p(p-1)}{2}}\Gamma(\alpha)\Gamma(\alpha-1)...\Gamma(\alpha-(p-1)),
\Re(\alpha)>p-1\eqno(1.1)
$$and
$$\tilde{\Gamma}_p(\alpha)=\int_{\tilde{X}>O}|{\rm det}(\tilde{X})|^{\alpha-p}{\rm e}^{-{\rm tr}(\tilde{X})}{\rm d}\tilde{X},\Re(\alpha)>p-1\eqno(1.2)
$$where $\Re(\cdot)$ means the real part of $(\cdot)$ and the integration is over all Hermitian positive definite matrix $\tilde{X}$. For our computations to follow, we will need some Jacobians of transformations in the complex domain. These will be listed here without proofs. For the proofs and for other such Jacobians, see Mathai (1997).

\vskip.3cm\noindent{\bf Lemma 1.1.\hskip.3cm}{\it Let $\tilde{X}$ and $\tilde{Y}$ be $m\times n$ with $mn$ distinct complex variables as elements. Let $A$ be $m\times m$ and $B$ be $n\times n$ nonsingular constant matrices. Then
$$\tilde{Y}=A\tilde{X}B,{\rm det}(A)\ne 0,{\rm det}(B)\ne 0\Rightarrow {\rm d}\tilde{Y}=[{\rm det}(A^{*}A)]^n[{\rm det}(B^{*}B)]^m{\rm d}\tilde{X}\eqno(1.3)
$$where $A^{*}$ and $B^{*}$ denote the conjugate transposes of $A$ and $B$ respectively; if $X,Y,A,B$ are real then
$$Y=AXB\Rightarrow {\rm d}Y=[{\rm det}(A)]^n[{\rm det}(B)]^m{\rm d}X\eqno(1.3a)
$$and if $a$ is a scalar quantity then
$$\tilde{Y}=a\tilde{X}\Rightarrow {\rm d}\tilde{Y}=|a|^{2mn}{\rm d}\tilde{X}.\eqno(1.3b)
$$}

\vskip.3cm\noindent{\bf Lemma 1.2.}\hskip.3cm{\it Let $\tilde{X}$ be $p\times p$ and Hermitian matrix of distinct complex variables as elements, except for Hermitianness. Let $A$ be a nonsingular constant matrix. Then
$$\tilde{Y}=A\tilde{X}A^{*}\Rightarrow{\rm d}\tilde{Y}=|{\rm det}(A)|^{-2p}{\rm d}\tilde{X}.\eqno(1.4)
$$If $A,X,Y,X=X'$ are real then
$$Y=AXA'\Rightarrow {\rm d}Y=[{\rm det}(A)]^{p+1}{\rm d}X.\eqno(1.4a)
$$If $Y,X,a,X=X'$ and $a$ scalar, then
$$Y=aX\rightarrow {\rm} dY=a^{\frac{p(p+1)}{2}}{\rm d}X\eqno(1.4c)
$$}

\vskip.3cm\noindent{\bf Lemma 1.3.}\hskip.3cm{\it Let $\tilde{X}$ be $p\times p$ and nonsingular with the regular inverse $\tilde{X}^{-1}$. Then
$$\tilde{Y}=\tilde{X}^{-1}\Rightarrow {\rm d}\tilde{Y}=\begin{cases}|{\rm det}(\tilde{X}^{*}\tilde{X})|^{-2p}{\rm d}\tilde{X}\mbox{ for a general $\tilde{X}$}\\
|{\rm det}(\tilde{X}^{*}\tilde{X})|^{-p}{\rm d}\tilde{X}\mbox{ for }\tilde{X}=\tilde{X}^{*}\mbox{ or }\tilde{X}=-\tilde{X}^{*}\end{cases}\eqno(1.5)
$$}
\vskip.3cm\noindent{\bf Lemma 1.4.}\hskip.3cm{\it Let $\tilde{X}$ be $p\times p$ Hermitian positive definite of distinct elements, except for Hermitian positive definiteness. Let $\tilde{T}$ be a lower triangular matrix where $\tilde{T}=(\tilde{t}_{ij}),\tilde{t}_{ij}=0,i<j, \tilde{t}_{ij},i\ge j$ are distinct, $\tilde{t}_{jj}=t_{jj}>0,j=1,...,p$, that is, the diagonal elements are real and positive. Then
$$\tilde{X}=\tilde{T}\tilde{T}^{*}\Rightarrow {\rm d}\tilde{X}=2^p\{\prod_{j=1}^pt_{jj}^{2(p-j)+1}\}{\rm d}\tilde{T}.\eqno(1.6)
$$}
\vskip.2cm With the help of Lemma 1.4 we can evaluate the complex matrix-variate gamma integral in (1.2) and show that it is equal to the expression in (1.1). When Lemma 1.4 is applied to the integral in (1.2) then the integral splits into $p$ integrals of the form
$$\prod_{j=1}^p2\int_0^{\infty}(t_{jj}^2)^{(\alpha-p)+\frac{1}{2}(2(p-j)+1)}{\rm e}^{-t_{jj}^2}{\rm d}t_{jj}=\prod_{j=1}^p\Gamma(\alpha-(j-1)),\Re(\alpha)>j-1,j=1,...,p
$$which results in the final condition as $\Re(\alpha)>p-1$, and $p(p-1)/2$ integrals of the form
\begin{align*}
\prod_{i>j}\int_{-\infty}^{\infty}{\rm e}^{-|\tilde{t}_{ij}|^2}{\rm d}\tilde{t}_{ij}
&=\prod_{i>j}\int_{-\infty}^{\infty}\int_{-\infty}^{\infty}{\rm e}^{-(t_{ij1}^2+t_{ij2}^2)}{\rm d}t_{ij1}\wedge{\rm d}t_{ij2}\\
&=\prod_{i>j}\sqrt{\pi}\sqrt{\pi}=\pi^{\frac{p(p-1)}{2}}, |\tilde{t}_{ij}|^2=t_{ij1}^2+t_{ij2}^2.
\end{align*}Thus the integral in (1.2) reduces to the expression in  (1.1).

\vskip.3cm\noindent{\bf Lemma 1.5.}\hskip.3cm{\it Let $\tilde{X}$ be $n\times p,n\ge p$ matrix of full rank $p$. Let $\tilde{S}=\tilde{X}^{*}\tilde{X}$, a $p\times p$ Hermitian positive definite matrix. Let ${\rm d}\tilde{X}$ and ${\rm d}\tilde{S}$ denote the wedge product of the differentials in $\tilde{X}$ and $\tilde{S}$ respectively. Then
$${\rm d}\tilde{X}=|{\rm det}(\tilde{S})|^{n-p}\frac{\pi^{np}}{\tilde{\Gamma}_p(n)}{\rm d}\tilde{S}.\eqno(1.7)
$$}
\vskip.2cm This is a very important result because $\tilde{X}$ is a rectangular matrix with $mn$ distinct elements whereas $\tilde{S}$ is Hermitian positive definite and $p\times p$. With the help of the above lemmas we will average a few functions over the Dirichlet measures in the complex domain.

\vskip.5cm{\section{Dirichlet Averages for Functions of Matrix Argument in the Complex Domain}}
\vskip.3cm

 All the matrices appearing in this section are $p\times p$ Hermitian positive definite unless stated otherwise. Consider the following complex matrix-variate type-1 Dirichlet measure:
\begin{align*}
f_1(\tilde{X}_1,...,\tilde{X}_k)&=\tilde{D}_k|{\rm det}(\tilde{X}_1)|^{\alpha_1-p}...|{\rm det}(\tilde{X}_k)|^{\alpha_k-p}\\
&\times |{\rm det}(I-\tilde{X}_1-...-\tilde{X}_k)|^{\alpha_{k+1}-p}\tag{2.1}
\end{align*}where $\tilde{X}_1,...\tilde{X}_k$ are $p\times p$ Hermitian positive definite, that is, $\tilde{X}_j>O,j=1,...,k$, such that $I-\tilde{X}_j>O,j=1,...,k,I-(\tilde{X}_1+...+\tilde{X}_k)>O$. The normalizing constant $\tilde{D}_k$ can be evaluated by integrating out matrices one at a time and the individual integrals are evaluated by using a complex matrix-variate type-1 beta integral of the form
$$\int_O^I|{\rm det}(\tilde{X})|^{\alpha-p}|{\rm det}(I-\tilde{X})|^{\beta-p}{\rm d}\tilde{X}=\frac{\tilde{\Gamma}_p(\alpha)\tilde{\Gamma}_p(\beta)}{\tilde{\Gamma}_p(\alpha+\beta)},
\Re(\alpha)>p-1,\Re(\beta)>p-1\eqno(2.2)
$$where $\tilde{\Gamma}_p(\alpha)$ is given in (1.1). It can be shown that the normalizing constant is the following:

$$\tilde{D}_k=\frac{\tilde{\Gamma}_p(\alpha_1+...+\alpha_{k+1})}{\tilde{\Gamma}_p(\alpha_1)...
\tilde{\Gamma}_p(\alpha_{k+1})}\eqno(2.3)
$$for $\Re(\alpha_j)>p-1,j=1,...,k+1$. Since (2.1), under (2.3) is a statistical density we can denote the averages of given functions as the expected values of those functions, denoted by $E(\cdot)$. Let  us consider a few functions and take their averages over the complex matrix-variate Dirichlet measure in (2.1). Let
$$\phi_1(\tilde{X}_1,...,\tilde{X}_k)=|{\rm det}(\tilde{X}_1)|^{\gamma_1}...|{\rm det}(\tilde{X}_k)|^{\gamma_k}.\eqno(2.4)
$$Then the average of (2.4) over the measure in (2.1) is given by
\begin{align*}
E[\phi_1]&=\tilde{D}_k\int_{\tilde{X}_1,...,\tilde{X}_k}|{\rm det}(\tilde{X}_1)|^{\alpha_1+\gamma_1-p}...|{\rm det}(\tilde{X}_k)|^{\alpha_k+\gamma_k-p}\\
&\times |{\rm det}(I-\tilde{X}_1-...-\tilde{X}_k)|^{\alpha_{k+1}-p}{\rm d}\tilde{X}_1\wedge....\wedge{\rm d}\tilde{X}_k.
\end{align*}Note that the only change is that $\alpha_j$ is changed to $\alpha_j+\gamma_j$ for $j=1,...,k$ and hence the result is available from the normalizing constant. That is,
$$E[\phi_1]=\{\prod_{j=1}^k\frac{\tilde{\Gamma}_p(\alpha_j+\gamma_j)}{\tilde{\Gamma}_p(\alpha_j)}\}
\frac{\tilde{\Gamma}_p(\alpha_1+...+\alpha_k)}
{\tilde{\Gamma}_p(\alpha_1+\gamma_1+...+\alpha_k+\gamma_k+\alpha_{k+1})},\eqno(2.5)
$$for $\Re(\alpha_j+\gamma_j)>p-1,j=1,...,k,\Re(\alpha_{k+1})>p-1$. Let
$$\phi_2(\tilde{X}_1,...,\tilde{X}_k)=|{\rm det}(I-\tilde{X}_1-...-\tilde{X}_k)|^{\delta}.\eqno(2.6)
$$Then in  the integral for $E[\phi_2]$ the only change is that the parameter $\alpha_{k+1}$ is changed to $\alpha_{k+1}+\delta$. Hence the result is available from the normalizing constant $\tilde{D}_k$. That is,
$$E[\phi_2]=\frac{\tilde{\Gamma}_p(\alpha_{k+1}+\delta)}{\tilde{\Gamma}_p(\alpha_{k+1})}
\frac{\tilde{\Gamma}_p(\alpha_1+...+\alpha_{k+1})}
{\tilde{\Gamma}_p(\alpha_1+...+\alpha_{k+1}+\delta)}\eqno(2.7)
$$for $\Re(\alpha_{k+1}+\delta)>p-1,\Re(\alpha_j)>p-1,j=1,...,k$. The structure in (2.7) is also the structure of the $\delta$-th moment of the determinant of the matrix having a complex matrix-variate type-1 beta distribution. Hence this $\phi_2$ has an equivalent representation in terms of the determinant of a matrix having a complex matrix-variate type-1 beta distribution.  Let
$$
\phi_3(\tilde{X}_1,...,\tilde{X}_k)={\rm e}^{{\rm tr}(\tilde{X}_1)}.\eqno(2.8)
$$Let us evaluate the Dirichlet average for $k=2$. Then

\begin{align*}
E[\phi_3]&=\tilde{D}_2\int_{\tilde{X}_1,\tilde{X}_2}{\rm e}^{{\rm tr}(\tilde{X}_1)}|{\rm det}(\tilde{X}_1)|^{\alpha_1-p}|{\rm det}(\tilde{X}_2)|^{\alpha_2-p}\\
&\times|{\rm det}(I-\tilde{X}_1-\tilde{X}_2)|^{\alpha_3-p}{\rm d}\tilde{X}_1\wedge...\wedge{\rm d}\tilde{X}_3.
\end{align*}
Take out $I-\tilde{X}_1$ from $|{\rm det}(I-\tilde{X}_1-\tilde{X}_2)|$ and make the transformation
$$\tilde{U}_2=(I-\tilde{X}_1)^{-\frac{1}{2}}\tilde{X}_2(I-\tilde{X}_1)^{-\frac{1}{2}}.
$$Then from Lemma 1.2, ${\rm d}\tilde{U}_2=|{\rm det}(I-\tilde{X}_1)|^{-p}{\rm d}\tilde{X}_2$. Now $\tilde{U}_2$ can be integrated out by using a complex matrix-variate type-1 beta integral given in (2.2). That is,
$$\int_{O<\tilde{U}_2<I}|{\rm det}(\tilde{U}_2)|^{\alpha_2-p}|{\rm det}(I-\tilde{U}_2)|^{\alpha_3-p}{\rm d}\tilde{U}_2=\frac{\tilde{\Gamma}_p(\alpha_2)\tilde{\Gamma}_p(\alpha_3)}{\tilde{\Gamma}_p(\alpha_2+\alpha_3)}\eqno(i)
$$for $\Re(\alpha_2)>p-1,\Re(\alpha_3)>p-1$. The $\tilde{X}_1$ integral to be evaluated is the following:
$$\int_{\tilde{X}_1}{\rm e}^{{\rm tr}(\tilde{X}_1)}|{\rm det}(\tilde{X}_1)|^{\alpha_1-p}|{\rm det}(I-\tilde{X}_1)|^{\alpha_2+\alpha_3-p}{\rm d}\tilde{X_1}.\eqno(ii)
$$In order to evaluate the integral in (ii) we can expand the exponential part by using zonal polynomials for complex argument, see Mathai (1997) and Mathai, Provost and Hayakawa (1995). We need a few notations and results from zonal polynomial expansions of determinants. The generalized Pochhammer symbol is the following:

$$[a]_M=\prod_{j=1}^p(a-j+1)_{k_j}=\frac{\tilde{\Gamma}_p(a,M)}{\tilde{\Gamma}_p(a)},\tilde{\Gamma}_p(a,M)
=\tilde{\Gamma}_p(a)[a]_M\eqno(2.9)
$$where the usual Pochhmmer symbol is
$$(a)_m=a(a+1)...(a+m-1),a\ne 0,(a)_0=1\eqno(2.10)
$$and $M$ represents the partition, $M=(m_1,...,m_p),m_1\ge m_2\ge ...\ge m_p,m_1+...+m_p=m$ and the zonal polynomial expansion for the exponential function is the following:
$${\rm e}^{{\rm tr}(\tilde{X})}=\sum_{m=0}^{\infty}\sum_M\frac{\tilde{C}_M(\tilde{X})}{m!}\eqno(2.11)
$$where $\tilde{C}_M(\tilde{X})$ is zonal polynomial of order $m$ in the complex matrix argument $\tilde{X}$, see (6.1.18) of Mathai (1997). One result on zonal polynomial that we require will be stated here as a lemma.

\vskip.3cm\noindent{\bf Lemma 2.1.}\hskip.3cm{\it
\begin{align*}
\int_{O<\tilde{Z}<I}&|{\rm det}(\tilde{Z})|^{\alpha-p}|{\rm det}(I-\tilde{Z})|^{\beta-p}\tilde{C}_M(\tilde{Z}\tilde{A}){\rm d}\tilde{Z}\\
&=\frac{\tilde{\Gamma}_p(\alpha,M)\tilde{\Gamma}_p(\beta)}{\tilde{\Gamma}_p(\alpha+\beta,M)}\tilde{C}_M(\tilde{A})\\
&=\frac{\tilde{\Gamma}_p(\alpha)\tilde{\Gamma}_p(\beta)}{\tilde{\Gamma}_p(\alpha+\beta)}\frac{(\alpha)_M}
{(\alpha+\beta)_M}
\tilde{C}_M(\tilde{A}),\tag{2.12}
\end{align*}}see also (6.1.21) of Mathai (1997), for $\Re(\alpha)>p-1,\Re(\beta)>p-1,\tilde{A}>O$. By using (2.12) we can evaluate the $
\tilde{X}_1$-integral in $E[\phi_3]$. That is,
\begin{align*}
\int_{O<\tilde{X}_1<I}&{\rm e}^{{\rm tr}(A\tilde{X}_1)}|{\rm det}(\tilde{X}_1)|^{\alpha_1-p}|{\rm det}(I-\tilde{X}_1)|^{\alpha_2+\alpha_3-p}{\rm d}\tilde{X}_1\\
&=\sum_{m=0}^{\infty}\sum_M\int_{O<\tilde{X}_1<I}\frac{\tilde{C}_M(\tilde{A}\tilde{X}_1)}{m!}|{\rm det}(\tilde{X}_1)|^{\alpha_1-p}|{\rm det}(I-\tilde{X}_1)|^{\alpha_2+\alpha_3-p}{\rm d}\tilde{X}_1\\
&=\sum_{m=0}^{\infty}\sum_M\frac{\tilde{C}_M(\tilde{A})}{m!}\frac{\tilde{\Gamma}_p(\alpha_1,M)
\tilde{\Gamma}_p(\alpha_2+\alpha_3)}
{\tilde{\Gamma}_p(\alpha_1+\alpha_2+\alpha_3,M)}.
\end{align*}Now, with the result on $\tilde{X}_2$-integral, $\tilde{D}_2$ and the above result will result in all the gamma products getting canceled and the final result is the following:
$$
E[\phi_3]=\sum_{m=0}^{\infty}\sum_M\frac{\tilde{C}_M(\tilde{A})}{m!}\frac{(\alpha_1)_M}
{(\alpha_1+\alpha_2+\alpha_3)_M}={_1F_1}(\alpha_1;\alpha_1+\alpha_2+\alpha_3;\tilde{A})\eqno(2.13)
$$for $\Re(\alpha_j)>p-1,j=1,2,3$ and ${_1F_1}$ is a confluent hypergeometric function of complex matrix argument $\tilde{A}$.

\vskip.5cm{\section{Dirichlet Averages in Complex Matrix-variate Type-2 Dirichlet Measure}}
\vskip.3cm
Consider the type-2 Dirichlet measure
\begin{align*}
f_2(\tilde{X}_1,...,\tilde{X}_k)&=\tilde{D}_k|{\rm det}(\tilde{X}_1)|^{\alpha_1-p}...|{\rm det}(\tilde{X}_k)|^{\alpha_k-p}\\
&\times |{\rm det}(I+\tilde{X}_1+...+\tilde{X}_k)|^{-(\alpha_1+...+\alpha_{k+1})}\tag{3.1}
\end{align*}for $\Re(\alpha_j)>p-1,j=1,...,k+1$ and it can be seen that the normalizing constant is the same as that in the type-1 Dirichlet measure. Let us evaluate some Dirichlet averages in the measure (3.1). Let
$$\phi_4(\tilde{X}_1,...,\tilde{X}_k)=|{\rm det}(\tilde{X}_1)|^{\gamma_1}...|{\rm det}(\tilde{X}_k)|^{\gamma_k}.\eqno(3.2)
$$Then when the average is taken the change is that $\alpha_j$ changes to $\alpha_j+\gamma_j,j=1,...,k$ and hence one should be table to find the value from the normalizing constant by adjusting for $\alpha_{k+1}$. Write $(\alpha_1+...+\alpha_{k+1})=(\alpha_1+\gamma_1+...+\alpha_k+\gamma_k)+(\alpha_{k+1}-\gamma_1-...-\gamma_k)$. That is, replace $\alpha_j$ by $\alpha_j+\gamma_j,j=1,...,k$ and replace $\alpha_{k+1}$ by $\alpha_{k+1}-\gamma_1-...-\gamma_k$ to obtain the result from the normalizing constant. Therefore
$$E[\phi_4]=\{\prod_{j=1}^k\frac{\tilde{\Gamma}_p(\alpha_j+\gamma_j)}{\tilde{\Gamma}_p(\alpha_j)}\}
\frac{\tilde{\Gamma}_p(\alpha_{k+1}-\gamma_1-...-\gamma_k)}{\tilde{\Gamma}_p(\alpha_{k+1})}\eqno(3.3)
$$for $\Re(\alpha_j+\gamma_j)>p-1,j=1,...,k$ and $\Re(\alpha_{k+1}-\gamma_1-...-\gamma_k)>p-1,\Re(\alpha_{k+1})>p-1$. Thus, only a few moments will exist, interpreting $E[\phi_4]$ as the product moment of the determinants of $\tilde{X}_1,...\tilde{X}_k$. Let
$$\phi_5(\tilde{X}_1,...,\tilde{X}_k)=|{\rm det}(I+\tilde{X}_1+...+\tilde{X}_k)|^{-\delta}.\eqno(3.4)
$$Then when the average is taken the only change in the integral is that $\alpha_{k+1}$ is changed to $\alpha_{k+1}+\delta$ and hence from the normalizing constant the result is the following:
$$E[\phi_5]=\frac{\tilde{\Gamma}_p(\alpha_{k+1}+\delta)}{\tilde{\Gamma}_p(\alpha_{k+1})}
\frac{\tilde{\Gamma}_p(\alpha_1+...+\alpha_{k+1})}
{\tilde{\Gamma}_p(\alpha_1+...+\alpha_{k+1}+\delta)},\eqno(3.5)
$$for $\Re(\alpha_{k+1}+\delta)>p-1$ and the other conditions on the parameters for $\tilde{D}_k$ remain the same. Observe that if $\Re(\delta)>0$ then the structure in (3.5) is that of the $\delta$-th moment of the determinant of a complex matrix-variate type-1 beta matrix. Thus, this type-2 form gives a type-1 form result. Let
$$\phi_6(\tilde{X}_1,\tilde{X}_2)={\rm e}^{-{\rm tr}(A\tilde{X}_1)}|{\rm det}(I+\tilde{X}_1)|^{\alpha_1+\alpha_3}.\eqno(3.6)
$$Then the Dirichlet average of $\phi_6$ in the complex matrix-variate type-2 Dirichlet measure in (3.1) for $k=2$ is the following:
\begin{align*}
E[\phi_6]&=\tilde{D}_2\int_{\tilde{X}_1,\tilde{X}_2}{\rm e}^{-{\rm tr}(\tilde{X}_1)}|{\rm det}(I+\tilde{X}_1)|^{\alpha_2+\alpha_3}|{\rm det}(\tilde{X}_1)|^{\alpha_1-p}|{\rm det}(\tilde{X}_2)|^{\alpha_2-p}\\
&\times |{\rm det}(I+\tilde{X}_1+\tilde{X}_2)|^{-(\alpha_1+\alpha_2+\alpha_3)}{\rm d}\tilde{X}_1\wedge...{\rm d}\tilde{X}_3.\end{align*}Take out $(I+\tilde{X}_1)$ from $I+\tilde{X}_1+\tilde{X}_2$ and make the transformation
$$\tilde{U}_2=(I+\tilde{X}_1)^{-\frac{1}{2}}\tilde{X}_2(I+\tilde{X}_1)^{-\frac{1}{2}}\Rightarrow{\rm d}\tilde{U}_2=|{\rm det}(I+\tilde{X}_1)|^{-p}{\rm d}\tilde{X}_2.
$$The $\tilde{U_2}$-integral gives
$$\int_{\tilde{U}_2>O}|{\rm det}(\tilde{U}_2)|^{\alpha_2-p}|{\rm det}(I+\tilde{U}_2)|^{-(\alpha_1+\alpha_2+\alpha_3)}{\rm d}\tilde{U}_2=\frac{\tilde{\Gamma}_p(\alpha_2)\tilde{\Gamma}_p(\alpha_1+\alpha_3)}
{\tilde{\Gamma}_p(\alpha_1+\alpha_2+\alpha_3)}.\eqno(i)
$$Observe that the exponent becomes zero and the factor containing $|{\rm det}(I+\tilde{X}_1)|$ disappears. Then the $\tilde{X}_1$-integral is
$$\int_{\tilde{X}_1>O}|{\rm det}(\tilde{X}_1)|^{\alpha_1-p}{\rm e}^{-{\rm tr}(A\tilde{X}_1)}{\rm d}\tilde{X}_1
=\tilde{\Gamma}_p(\alpha_1)|{\rm det}(A)|^{-\alpha_1}.\eqno(ii)
$$The results from (i), (ii) and $\tilde{D}_2$ gives the final result as follows:
$$E[\phi_6]=\frac{\tilde{\Gamma}_p(\alpha_1+\alpha_3)}{\tilde{\Gamma}_p(\alpha_3)}|{\rm det}(A)|^{-\alpha_1}\eqno(3.7)
$$and the original conditions on the parameters remain the same and no further conditions are needed, where $A>O$. Note that if $\phi_6$ did not have the factor $|{\rm det}(I+\tilde{X}_1)|^{\alpha_1+\alpha_3}$ then a factor containing $|{\rm det}(I+\tilde{X}_1)|$ would also have been present then the $\tilde{X}_1$-integral would have gone in terms of a Whittaker function of matrix argument, see Mathai (1997).

\vskip.5cm{\section{Dirichlet Averages in Complex Rectangular Matrix-variate Dirichlet Measure}}
\vskip.3cm

Let $B_j$ be $n_j\times n_j$ Hermitian positive definite constant matrix and let $B_j^{\frac{1}{2}}$ denote the Hermitian positive definite square root of $B_j$. Let $\tilde{X}_j$ be $n_j\times p$, $n_j\ge p$ matrix of full rank $p$ so that $\tilde{X}_j^{*}\tilde{X}_j=\tilde{S}_j>O$ or $\tilde{S}_j$ is Hermitian positive definite. Observe that for $p=1$, $\tilde{X}_j^{*}B_j\tilde{X}_j$ is a positive definite Hermitian form. Hence our results to follow will also cover results on Hermitian forms. Consider the model
\begin{align*}
f_3(\tilde{X}_1,...,\tilde{X}_k)&=\tilde{G}_k|{\rm det}(\tilde{X}_1^{*}B_1\tilde{X}_1)|^{\alpha_1}...|{\rm det}(\tilde{X}_k^{*}B_k\tilde{X}_k)|^{\alpha_k}\\
&\times |{\rm det}(I-\tilde{X}_1^{*}B_1\tilde{X}_1-...-\tilde{X}_k^{*}B_k\tilde{X}_k)|^{\alpha_{k+1}-p}\tag{4.1}
\end{align*}where $\tilde{G}_k$ is the normalizing constant and $O<\tilde{X}_j^{*}B_j\tilde{X}_j<I,j=1,...,k,O<\tilde{X}_1^{*}B_1\tilde{X}_1+...+\tilde{X}_k^{*}B_k\tilde{X}_k<I,
j=1,...,k$. The normalizing constant is evaluated by using the following procedure. Let $\tilde{Y}_j=B_j^{\frac{1}{2}}\tilde{X}_j\Rightarrow {\rm d}\tilde{Y_j}=|{\rm det}(B_j)|^p{\rm d}\tilde{X}_j$ from Lemma 1.1. Let $\tilde{Y}_j^{*}\tilde{Y}_j=\tilde{S}_j$. Then from Lemma 1.5 we have
$${\rm d}\tilde{Y}_j=\frac{\pi^{n_jp}}{\tilde{\Gamma}_p(n_j)}|{\rm det}(\tilde{S}_j)|^{n_j-p}{\rm d}\tilde{X}_j.\eqno(i)
$$Then
$${\rm d}\tilde{X}_1\wedge...\wedge{\rm d}\tilde{X}_k=\{\prod_{j=1}^k\frac{\pi^{n_jp}}{\tilde{\Gamma}_p(n_j)}|{\rm det}(\tilde{B}_j)|^{-p}|{\rm det}(\tilde{S}_j)|^{n_j-p}\}{\rm d}\tilde{S}_1\wedge...\wedge{\rm d}\tilde{S}_k.\eqno(ii)
$$Since the total integral is $1$ we have
\begin{align*}
1&=\int_{\tilde{X}_1,...\tilde{X}_k}f_3(\tilde{X}_1,...,\tilde{X}_k){\rm d}\tilde{X}_1\wedge...\wedge{\rm d}\tilde{X}_k\\
&=\tilde{G}_k\{\prod_{j=1}^k\frac{\pi^{n_jp}}{\tilde{\Gamma}_p(n_j)}\}\int_{\tilde{S}_1,...,\tilde{S}_k}|{\rm det}(\tilde{S}_1)|^{\alpha_1+n_1-p}...\\
&\times |{\rm det}(\tilde{S}_k)|^{\alpha_k+n_k-p}|{\rm det}(I-\tilde{S}_1-...-\tilde{S}_k)|^{\alpha_{k+1}-p}{\rm d}\tilde{S}_1\wedge...\wedge{\rm d}\tilde{S}_k.
\end{align*}Now, evaluating the type-1 Dirichlet integrals over the $\tilde{S}_j$'s one has the result. That is,
\begin{align*}
\tilde{G}_k&=\{\prod_{j=1}^k|{\rm det}(B_j)|^{p}\frac{\tilde{\Gamma}_p(n_j)}{\pi^{n_jp}}\frac{1}{\tilde{\Gamma}_p(\alpha_j+n_j)}\}\\
&\times \frac{\tilde{\Gamma}_p(\alpha_1+...+\alpha_{k+1}+n_1+...+n_k)}{\tilde{\Gamma}_p(\alpha_{k+1})}\tag{4.2}\end{align*}
for $B_j>O,\Re(\alpha_j+n_j)>p-1,j=1,...,k,\Re(\alpha_{k+1})>p-1$. Thus, (4.1) with (4.2) defines a rectangular complex matrix-variate type-1 Dirichlet measure. Thee is a corresponding type-2 Dirichlet measure, given by the following:
\begin{align*}
f_4(\tilde{X}_1,...,\tilde{X}_k)&=\tilde{G}_k|{\rm det}(\tilde{X}_1)|^{\alpha_1}...|{\rm det}(\tilde{X}_k)|^{\alpha_k}\\
&\times |{\rm det}(I+\tilde{X}_1+...+\tilde{X}_k)|^{-(\alpha_1+...+\alpha_{k+1}+n_1+...+n_k)}\tag{4.3}
\end{align*}for $B_j>O,\Re(\alpha_j+n_j)>p-1,j=1,...,k,\Re(\alpha_{k+1})>p-1$ and $\tilde{G}_k$ is the same as the one appearing in (4.2). Let us compute the Dirichlet averages of some functions in the type-2 rectangular complex matrix-variate Dirichlet measure in (4.3). Let
$$\phi_7(\tilde{X}_1,...,\tilde{X}_k)=|{\rm det}(\tilde{X}_1)|^{\gamma_1}...|{\rm det}(\tilde{X}_k)|^{\gamma_k}.\eqno(4.4)
$$Then when we take the expected value of $\phi_7$ in (4.3) the only change is that $\alpha_j$ changes to $\alpha_j+\gamma_j,j=1,...,k$ and hence the final result is available from the normalizing constant. Therefore
$$E[\phi_7]=\{\prod_{j=1}^k\frac{\tilde{\Gamma}_p(\alpha_j+n_j+\gamma_j)}{\tilde{\Gamma}_p(\alpha_j+n_j)}\}
\frac{\tilde{\Gamma}_p(\alpha_{k+1}-\gamma_1-...-\gamma_k)}{\tilde{\Gamma}_p(\alpha_{k+1})}\eqno(4.5)
$$for $\Re(\alpha_j+n_j+\gamma_j)>p-1,j=1,...,k,\Re(\alpha_{k+1}-\gamma_1-...-\gamma_k)>p-1,\Re(\alpha_{k+1})>p-1$.
Let
$$\phi_8(\tilde{X}_1,...,\tilde{X}_k)=|{\rm det}(I+\tilde{X}_1+...+\tilde{X}_k)|^{-\delta}.\eqno(4.6)
$$Then the only change is that $\alpha_{k+1}$ goes to $\alpha_{k+1}+\delta$ in the integral and no other change is there and hence the average is available from the normalizing constant. That is,
$$E[\phi_8]=\frac{\tilde{\Gamma}_p(\alpha_{k+1}+\delta)}{\tilde{\Gamma}_p(\alpha_{k+1})}
\frac{\tilde{\Gamma}_p(\alpha_1+...+\alpha_{k+1}+n_1+...+n_k)}
{\tilde{\Gamma}_p(\alpha_1+...+\alpha_{k+1}+n_1+...+n_k+\delta)}\eqno(4.7)
$$for $\Re(\alpha_j+n_j)>p-1,j=1,...,k,\Re(\alpha_{k+1}+\delta)>p-1,\Re(\alpha_{k+1})>p-1$.
\vskip.2cm The case $p=1$ in the complex rectangular matrix-variate type-1 Dirichlet measure is very interesting. We have a set of Hermitian positive definite quadratic forms here having a joint density of the following form:

\begin{align*}
f_5(\tilde{X}_1,...,\tilde{X}_k)&=\tilde{G}_k[\tilde{X}_1^{*}B_1\tilde{X}_1]^{\alpha_1}...
[\tilde{X}_k^{*}B_k\tilde{X}_k]^{\alpha_k}\\
&\times |{\rm det}(I-[\tilde{X}_1^{*}B_1\tilde{X}_1]-...-[\tilde{X}_k^{*}B_k\tilde{X}_k])|^{\alpha_{k+1}-p}\tag{4.8}
\end{align*}where $B_j>O$, and $\tilde{X}_j^{*}B_j\tilde{X}_j$ is a scalar quantity, $j=1,...,k$. Consider the same types of transformations as before. $\tilde{Y}_j=B_j^{\frac{1}{2}}\tilde{X}_j$. Then $\tilde{Y}_j^{*}\tilde{Y}_j=|\tilde{y}_{j1}|^2+...+|\tilde{y}_{jn_j}|^2$ or the sum or squares of the absolute values of $\tilde{y}_{jr}$ where $\tilde{Y}_j^{*}=(\tilde{y}_{j1}^{*},...,\tilde{y}_{jn_j}^{*})$. This is an isotropic point in in the $2n_j$-dimensional Euclidean space. From here, one can establish various connections to geometrical probability problems, see Mathai (1999). Also (4.8) is associated with the theory of generalized Hermitian forms in pathway models, see Mathai (2007). Let us evaluate the $h$-th moment of
$$\phi_9(\tilde{X}_1,...,\tilde{X}_k)=[\tilde{X}_1^{*}B_1\tilde{X}_1+...+\tilde{X}_k^{*}B_k\tilde{X}_k]^h\eqno(4.9)
$$for $p=1$. For $p>1$ we have seen that this is not available directly but moments of $|{\rm det}(I-\tilde{X}_1^{*}B_1\tilde{X}_1-...-\tilde{X}_k^{*}B_k\tilde{X}_k)|$ was available. But for $p=1$ one can obtain $h$-th moment of both for an arbitrary $h$. By computing the $h$-th moment of $[1-\tilde{X}_1^{*}B_1\tilde{X}_1-...-\tilde{X}_k^{*}B_k\tilde{X}_k]$ for $p=1$ we note that  for arbitrary $h$  this quantity and its complementary part $[\tilde{X}_1^{*}B_1\tilde{X}_1+...+\tilde{X}_k^{*}B_k\tilde{X}_k]$ are both  scalar variable type-1 beta distributed with the parameters $(\alpha_{k+1},\sum_{j=1}^k(\alpha_j+n_j))$ and $(\sum_{j=1}^k(\alpha_j+n_j),\alpha_{k+1})$ respectively. Then
$$E[\phi_9]=\frac{\tilde{\Gamma}_p(\sum_{j=1}^k(\alpha_j+n_j)+h)}{\tilde{\Gamma}_p(\sum_{j=1}^k(\alpha_j+n_j))}
\frac{\tilde{\Gamma}_p(\sum_{j=1}^k(\alpha_j+n_j)+\alpha_{k+1})}{\tilde{\Gamma}_p(\sum_{j=1}^k(\alpha_j+n_j)+\alpha_{k+1}+h)}\eqno(4.10)
$$for $\Re(\alpha_j)>p-1,j=1,...,k+1,\Re(\sum_{j=1}^k(\alpha_j+n_j)+h)>p-1$. Consider $\phi_9$ in the complex matrix-variate type-2 Dirichlet measure for $p=1$. Then the $h$-th moment will reduce to the following:
$$E[\phi_9]=\frac{\tilde{\Gamma}_p(\sum_{j=1}^k(\alpha_j+n_j)+h)}{\tilde{\Gamma}_p(\sum_{j=1}^k(\alpha_j+n_j))}
\frac{\tilde{\Gamma}_p(\alpha_{k+1}-h)}{\tilde{\Gamma}_p(\alpha_{k+1})}\eqno(4.11)
$$for $\Re(\alpha_{k+1}-h)>p-1,\Re(\alpha_j)>p-1,j=1,...,k+1,\Re(\sum_{j=1}^k(\alpha_j+n_j)+h)>p-1$.
\vskip.2cm Many such results can be obtained for the type-1 and type-2 Dirichlet measures in Hermitian positive definite Dirichlet measures or in rectangular matrix-variate Dirichlet measures.

\vskip.5cm{\section{Applications}}
\vskip.3cm
For our applications in the theory of special functions, fractional calculus, biology, probability, and stochastic processes, Dirichlet averages and their diverse approaches are used.  In this section, the main areas where the applications of
Dirichlet averages are presented:
\subsection{Special Functions}
Dirichlet averages were introduced by Carlson in his work Carlson (1977).Carlson (1963,1969, 1975, 1987) observed that the straightforward idea of this kind of averaging generalizes and unifies a wide range of special functions, including various orthogonal polynomials and generalized hyper-geometric functions. The relationship between Dirichlet splines and an important
class of hypergeometric functions of several variables is given in Neuman and Fleet (1994), and Carlson (1991).
Numerous investigations of B-splines, including those by Carlson (1991), Massopust and Forster (2010), and Stolarsky means, by Simi\'{c} and Bin-Mohsin (2020) used Dirichlet averages.

\subsection{Fractional Calculus}
The
Dirichlet average of elementary functions like power
function, exponential function, etc. is given by many
notable mathematicians. There are many results available in the literature converting the elementary function into the summation form after
taking the Dirichlet average of those functions, using fractional
integral, and getting new results, see Kilbas and Kattuveettill (2008), Saxena et al.(2010), Kumar et al. (2022). Those results will be used in
the future by mathematicians and scientists in a variety of fields.
\subsection{Gene Expression Modelling}
Clustering is a key data processing technique for interpreting microarray data and determining genetic networks. Hierarchical Dirichlet processes (HDP) clustering is able to capture the hierarchical elements that are common in biological data, such as gene expression data, by including a hierarchical structure into the statistical model. Wang and Wang (2013) presented a hierarchical Dirichlet process model for gene expression clustering.
\subsection{Geometrical Probability}
Thomas and Mathai (2009) propose a generalized Dirichlet model application to geometrical probability problems.
When the linearly independent random points in Euclidean $n$ space have highly general real
rectangle matrix-variate beta density, the volumes of random parallelotopes are explored.
In order to evaluate statistical hypotheses,
structural decomposition is provided, and random volumes are linked to generalized Dirichlet models and likelihood ratio criteria.
This makes it possible to calculate percentage points of random volumes using the generalized Dirichlet marginal's $p-$values.
\subsection{Bayesian Analysis}
Carlson's original definition of Dirichlet averages is expressed as mixed multinomial distributions' probability-generating functions.They also contribute significantly to the solution of elliptic integrals and have several connections to statistical applications. Dickey (1983) obtained that several nested families are built for Bayesian inference in multinomial sampling and contingency tables that generalize the Dirichlet distributions.  These distributions can be used to model populations of personal probabilities evolving under the process of inference from statistical data.


\end{document}